\documentclass[11pt,a4paper]{amsart}
\usepackage{amssymb}

\DeclareMathOperator{\PSL}{PSL}
\DeclareMathOperator{\PGL}{PGL}
\DeclareMathOperator{\SL}{SL}
\DeclareMathOperator{\GL}{GL}
\DeclareMathOperator{\Sz}{Sz}

\theoremstyle{plain}
\newtheorem{theorem}{Theorem}[section]
\newtheorem{prop}[theorem]{Proposition}
\newtheorem{lemma}[theorem]{Lemma}
\newtheorem{corr}[theorem]{Corollary}

\theoremstyle{definition}
\newtheorem{Theorem}{Theorem}

\newtheorem{conj}{Conjecture}

\numberwithin{equation}{section}

\title[Connectivity of the PRA Graph of $\PSL(2,q)$]
{Connectivity of the Product Replacement Algorithm Graph of
$\bf{\PSL(2,q)}$}
\author{Shelly Garion}
\address{Einstein Institute of Mathematics,
The Hebrew University of Jerusalem, Jerusalem 91904, Israel}
\email{shellyg@math.huji.ac.il}
\begin{document}

%%%%%%%%%%%%%%%%%%%%%%%%%%%%%%%%%%%%%%%%%%%%%%%%%%%%%%%

\begin{abstract}
The product replacement algorithm is a practical algorithm to
construct random elements of a finite group $G$. It can be described
as a random walk on a graph $\Gamma_k(G)$ whose vertices are the
generating $k$-tuples of $G$ (for a fixed $k$). We show that if
$G=\PSL(2,q)$ or $\PGL(2,q)$, where $q$ is a prime power, then
$\Gamma_k(G)$ is connected for any $k \geq 4$. This generalizes
former results obtained by Gilman and Evans.
\end{abstract}

%\footnotetext{2000 {\it Mathematics Subject Classification:} 20D06,
%20D60. }


\footnotetext{This article was submitted to the {\it Journal of
Group Theory} on July 2007 and accepted on December 2007.}

\maketitle

%%%%%%%%%%%%%%%%%%%%%%%%%%%%%%%%%%%%%%%%%%%%%%%%%%%%%%%

\section{Introduction}

%%%%%%%%%%%%%%%%%%%%%%%%%%%%%%%%%%%%%%%%%%%%%%%%%%%%%%%

\subsection{The product replacement algorithm}
The \emph{product replacement algorithm} (PRA) is a practical
algorithm to construct random elements of a finite group. The
algorithm was introduced and analyzed in~\cite{CLMNO}, where the
authors proved that it produces asymptotically uniformly distributed
elements. They also showed that the algorithm has very good
performance in several interesting cases. As the success of the
algorithm has become widely acknowledged, it has been included as a
standard routine in the two major algebra packages \textsf{GAP} and
\textsc{Magma}. Since then the algorithm has been widely
investigated (see~\cite{BP,GaP,LP,P}).

The product replacement algorithm is defined as follows. Let $G$ be
a finite group and let $d(G)$ be the minimal number of generators of
$G$. For any integer $k \geq d(G)$, let
\[
V_k(G)=\{(g_1,\ldots,g_k) \in G^k: \langle g_1,\ldots,g_k \rangle =
G\}
\]
be the set of all \emph{generating $k$-tuples} of $G$. Given a
generating $k$-tuple, a \emph{move} to another such tuple is defined
by first uniformly selecting a pair $(i,j)$ with $1 \leq i \neq j
\leq k$ and then applying one of the following four operations with
equal probability:
\[
\begin{split}
    R_{i,j}^{\pm}:  (g_1,\ldots,g_i,\ldots,g_k) \mapsto
    (g_1,\ldots,g_i\cdot g_j^{\pm 1},\ldots,g_k) \\
    L_{i,j}^{\pm}:  (g_1,\ldots,g_i,\ldots,g_k) \mapsto
    (g_1,\ldots,g_j^{\pm 1} \cdot g_i,\ldots,g_k)
\end{split}
\]

To produce a random element in $G$, start with some generating
$k$-tuple, apply the above moves several times, and finally return a
random element of the generating $k$-tuple that was reached. In
practice, one applies approximately 100 moves.

%%%%%%%%%%%%%%%%%%%%%%%%%%%%%%%%%%%%%%%%%%%%%%%%%%%%%%%

\subsection{The product replacement algorithm graph}
The moves in the PRA can be conveniently encoded by the \emph{PRA
graph} $\Gamma_k(G)$ whose vertices are the tuples $V_k(G)$, with
edges corresponding to the moves $R_{i,j}^{\pm}, L_{i,j}^{\pm}$. The
PRA corresponds to a random walk on this graph. We are interested in
studying the connectivity of this graph.

However, it is usually more convenient to look at the \emph{extended
PRA graph} $\tilde \Gamma_k(G)$. This is a graph on $V_k(G)$
corresponding to the so called \emph{Nielsen moves:} $R_{i,j}^{\pm},
L_{i,j}^{\pm}$ and $P_{i,j}, I_{i}$, $1 \leq i \neq j \leq k$, where
\[
\begin{split}
    P_{i,j}&:  (g_1,\ldots,g_i,\ldots,g_j,\ldots,g_k) \mapsto
    (g_1,\ldots,g_j,\ldots,g_i,\ldots,g_k)  \\
    I_{i}&:  (g_1,\ldots,g_i,\ldots,g_k) \mapsto
    (g_1,\ldots,g_i^{-1},\ldots,g_k)
\end{split}
\]

It is clear from the definitions that if  $\chi_k(G)$ and
$\tilde\chi_k(G)$ denote the number of connected components in
$\Gamma_k(G)$ and $\tilde \Gamma_k(G)$ respectively, then $\tilde
\chi_k(G) \leq \chi_k(G) \leq 2\tilde \chi_k(G) $. Moreover, if
$k\geq d(G)+1$ then $\Gamma_k(G)$ is connected if and only if
$\tilde \Gamma_k(G)$ is connected \cite[Prop. 2.2.1]{P}.

It is worthwhile to mention here the relation to the notion of
\emph{$T$-systems}. Originally, $T$-systems were introduced in
\cite{NN} to study presentations of finite groups. However, it
became apparent that many results that were originally obtained for
$T$-systems can be restated for the extended PRA graph as well (see
\cite{Du1,Du2,E,Gi}, for example).

%%%%%%%%%%%%%%%%%%%%%%%%%%%%%%%%%%%%%%%%%%%%%%%%%%%%%%%

\subsection{Connectivity of the PRA graph}

Let $G$ be a finite group and let $k \geq d(G)$ be an integer. In
this paper we consider the connectivity of the extended PRA graph
$\tilde \Gamma_k(G)$.

There are several examples for which $k=d(G)$ and $\tilde
\Gamma_k(G)$ is not connected (see~\cite{Du1,GS,GP,Ne,P}).
Pak~\cite{P} asked whether there is a finite group $G$ and an
integer $k\geq d(G)+1$ such that $\tilde \Gamma_k(G)$ is
disconnected. As there are no known such examples, the following
conjecture is naturally raised, see~\cite{CP}.

\begin{conj} \label{conj_gen}
If $k\geq d(G)+1$ then $\tilde \Gamma_k(G)$ is connected.
\end{conj}

This conjecture was first proved in~\cite{NN} for finite abelian
groups, and later it was generalized by Dunwoody~\cite{Du2} to
finite solvable groups. Gilman \cite[Thm. 3]{Gi} proved that for any
finite group $G$, $\tilde \Gamma_k(G)$ is connected if $k \geq 2
\log_2(|G|)$. Subsequently, a better bound has been obtained by
Pak~\cite{P}: if $\mu(G)$ is the maximal size of an
\emph{independent generating set} of $G$, i.e. a set of generators
such that no member can be omitted, then $\tilde \Gamma_k(G)$ is
connected for any $k \geq d(G)+\mu(G)$.

%%%%%%%%%%%%%%%%%%%%%%%%%%%%%%%%%%%%%%%%%%%%%%%%%%%%%%%

\subsection{Simple groups}

It is well known that $d(G)=2$ for any non-abelian finite simple
group $G$. In this special case Conjecture~\ref{conj_gen} becomes

\begin{conj} \textbf{(Wiegold).}
If $G$ is a finite simple group and $k\geq 3$, then $\tilde
\Gamma_k(G)$ is connected.
\end{conj}

This conjecture has been proved only in the following cases.
\begin{prop} \label{simple_conn}
$\tilde \Gamma_k(G)$ is connected in the following cases:
\begin{enumerate}
\renewcommand{\theenumi}{\alph{enumi}}
\item \label{t:PSLp} %\cite[Thm. 1]{Gi}.
$G=\PSL(2,p)$, where $p\geq 5$ is prime and $k \geq 3$.
\item \label{t:PSL2} %\cite{E}.
$G=\PSL(2,2^m)$, where $m\geq 2$ and $k \geq 3$.
\item \label{t:Su} %\cite{E}.
$G=\Sz(2^{2m-1})$, where $m\geq 2$ and $k \geq 3$.
\item \label{t:An} %\cite{CP, Da}.
$G=A_n$, where $6 \leq n \leq 11$ and $k=3$.
\end{enumerate}
\end{prop}

\begin{proof}
Part~(\ref{t:PSLp}) is~\cite[Thm. 1]{Gi}. See~\cite{E}
for~(\ref{t:PSL2}) and~(\ref{t:Su}), while~(\ref{t:An}) follows from
results in~\cite{CP} and~\cite{Da}. Note that the proof in~\cite{CP}
is based on computer calculations.
\end{proof}

The aim of this paper is to extend~(\ref{t:PSLp}) and
~(\ref{t:PSL2}) in Proposition~\ref{simple_conn} to $\PSL(2,q)$ for
$q=p^e$, where $p$ is an odd prime and $e > 1$.

A naive bound can be easily computed in view of the following
observation. By~\cite{WS}, $\mu(\PSL(2,q)) \leq \max(6, \pi+2)$,
where $\pi = \pi(e)$ is the number of distinct prime divisors of
$e$, hence the aforementioned result of Pak~\cite{P} implies that
$\tilde \Gamma_k(\PSL(2,q))$ is connected for $k \geq \max(8,
\pi+4)$. However, this bound still depends on $q$. In the following
theorem we present a bound which is independent of $q$.

\begin{Theorem} \label{main}
Let $G=\PSL(2,q)$ or $\PGL(2,q)$, where $q$ is an odd prime power.
Then $\tilde \Gamma_k(G)$ is connected for any $k \geq 4$.
\end{Theorem}

In Section 2, we present some of the basic properties of the groups
$\PSL(2,q)$ and $\PGL(2,q)$ that are needed in the proof of
Theorem~\ref{main}. The proof itself is presented in Section 3.
Here, we adapt some of the techniques of Evans~\cite{E} to the
groups $\PSL(2,q)$ over arbitrary finite fields.

%%%%%%%%%%%%%%%%%%%%%%%%%%%%%%%%%%%%%%%%%%%%%%%%%%%%%%%

\medskip

\emph{Acknowledgements.} This paper is part of the author's Ph.D.
studies under the guidance of Alex Lubotzky, whom I would like to
thank for his support and good advice. I am grateful to Avinoam Mann
and Aner Shalev for many useful discussions. I would also like to
thank the referee for the many helpful comments and suggestions.

The author acknowledges the support of the Israeli Ministry of
Science, Culture and Sport.

%%%%%%%%%%%%%%%%%%%%%%%%%%%%%%%%%%%%%%%%%%%%%%%%%%%%%%%

\section{Preliminaries}

\subsection{Definitions}
Let $q=p^e$, where $p$ is an odd prime and $e \geq 1$. Recall that
$\GL(2,q)$ is the group of invertible $2 \times 2$ matrices over the
finite field with $q$ elements, which we denote by $\mathbb{F}_q$,
and $\SL(2,q)$ is the subgroup of $\GL(2,q)$ comprising the matrices
with determinant $1$. Then $\PGL(2,q)$ and $\PSL(2,q)$ are the
quotients of $\GL(2,q)$ and $\SL(2,q)$ by their respective centers.
The orders of $\PGL(2,q)$ and $\PSL(2,q)$ are $q(q-1)(q+1)$ and
$\frac{1}{2}q(q-1)(q+1)$ respectively, therefore we can identify
$\PSL(2,q)$ with a normal subgroup of index $2$ in $\PGL(2,q)$. Also
recall that $\PSL(2,q)$ is simple for $q \neq 3$.

Let $\mathbb{P}_1(q)$ denote the projective line over
$\mathbb{F}_q$. Then $\PGL(2,q)$ acts on $\mathbb{P}_1(q)$ by
\[
    \begin{pmatrix}
    a & b \\ c&d
    \end{pmatrix}: \quad
    z \mapsto \frac{az+b}{cz+d}
\]
hence, it can be identified with the group of projective
transformations on $\mathbb{P}_1(q)$. Under this identification,
$\PSL(2,q)$ is the set of all transformations for which $ad-bc$ is a
square in $\mathbb{F}_q$. It is well known that $\PGL(2,q)$ is
triply transitive on $\mathbb{P}_1(q)$, while $\PSL(2,q)$ acts
2-transitively.

%%%%%%%%%%%%%%%%%%%%%%%%%%%%%%%%%%%%%%%%%%%%%%%%%%%%%%%

\subsection{Group elements} \label{elements}

One can classify the elements of $\PGL(2,q)$ and $\PSL(2,q)$
according to their action on $\mathbb{P}_1(q)$. This is the same as
considering the possible Jordan forms. The following table lists the
three types of elements according to whether they have $0,1$ or $2$
distinct eigenvalues in $\mathbb{F}_q$.

\begin{table}[h]
\begin{tabular} {|c|c|c|c|}
\hline
type & action on $\mathbb{P}_1(q)$ & order in $\PGL(2,q)$ & order in $\PSL(2,q)$  \\
\hline

unipotent & fixes $1$ point & $p$ & $p$ \\

split & fixes $2$ points & divides $q-1$ & divides $\frac{1}{2}(q-1)$  \\

non-split & no fixed points & divides $q+1$ & divides $\frac{1}{2}(q+1)$ \\
\hline
\end{tabular}
\end{table}

%%%%%%%%%%%%%%%%%%%%%%%%%%%%%%%%%%%%%%%%%%%%%%%%%%%%%%%

\subsection{Subgroups} \label{subgroups}

The classification of subgroups of $\PSL(2,q)$ and $\PGL(2,q)$ is
well known, and is originally due to Dickson~\cite{D}
(cf.~\cite{Hu,Su}). The following table specifies all the subgroups
of $\PSL(2,q)$ and $\PGL(2,q)$ up to isomorphism, divided into the
various Aschbacher classes, following~\cite{As}. Here $w_1$ and
$w_2$ are split and non-split elements of $G$, respectively.

\begin{table}[h]
\begin{tabular} {|c|c|c|c|c|}
\hline
 &  & maximal & maximal &
 \\
class & type &  order in & order in & structure and properties \\
 & &  $\PGL(2,q)$ & $\PSL(2,q)$ & \\
\hline \hline

$\mathcal{C}_1$ & $p$-group & $q$ & $q$ & elementary abelian \\
& & & & $p$-group \\
\cline{2-5}

 & Borel & $q(q-1)$ & $\frac{1}{2}q(q-1)$ & normalizer of a
$p$-group
\\
& & & & stabilizes a point in $\mathbb{P}_1(q)$ \\
\hline

$\mathcal{C}_2$ & cyclic & $q-1$ & $\frac{1}{2}(q-1)$ & $\langle w_1 \rangle$ \\
\cline{2-5}

& dihedral & $2(q-1)$ & $q-1$ & normalizer of $\langle w_1 \rangle$ \\
\hline

$\mathcal{C}_3$ & cyclic & $q+1$ & $\frac{1}{2}(q+1)$ & $\langle w_2 \rangle$ \\
\cline{2-5}

& dihedral & $2(q+1)$ & $q+1$ & normalizer of $\langle w_2 \rangle$ \\
\hline

$\mathcal{C}_5$ & $\PSL(2,q_1)$ & -- & --  & exists if $q=q_1^r$ $(r \in \mathbb{N})$\\
\cline{2-5}

& $\PGL(2,q_1)$ & -- & --  & exists in $\PGL(2,q)$ \\
& & & & if $q=q_1^r$ $(r \in \mathbb{N})$\\
& & & & and in $\PSL(2,q)$ \\
& & & & if $q=q_1^{2r}$ $(r \in \mathbb{N})$\\
\hline

$\mathcal{C}_6$ & $A_4$ & 12 & 12 & -- \\
\cline{2-5}

& $S_4$ & 24 & 24 & exists in $\PSL(2,q)$ if \\
& & & & $q \equiv \pm 1 \pmod 8$ \\
\hline

$\mathcal{S}$ & $A_5$ & 60 & 60 & exists if $p=5$ \\
& & & & or $q \equiv \pm 1 \pmod {10}$ \\
\hline

\end{tabular}
\end{table}

Since the maximal subgroups in the classes $\mathcal{C}_1$,
$\mathcal{C}_2$ and $\mathcal{C}_3$ arise from stabilizers of
subspaces of $\mathbb{F}_q$ or $\mathbb{F}_q^2$, we will call them
\emph{structural} subgroups. The subgroups of class $\mathcal{C}_5$
are usually called \emph{subfield} subgroups. For convenience we
will refer to the subgroups $A_4$, $S_4$ and $A_5$ as \emph{small}.

In order to study the action of these subgroups on
$\mathbb{P}_1(q)$, we introduce the following convenient notation.
Let $G$ be a group acting on $\mathbb{P}_1(q)$ and let $a,b \in
\mathbb{P}_1(q)$ be two distinct points. Denote by $G_a$ the
subgroup of $G$ that fixes $a$, and by $G_{a,b}$ the subgroup of $G$
that fixes $a$ and $b$ pointwise.

%%%%%%%%%%%%%%%%%%%%%%%%%%%%%%%%%%%%%%%%%%%%%%%%%%%%%%%

\subsection{Conjugation of elements and subgroups}

In~\cite{D}, Dickson actually classifies all the subgroups of
$\PSL(2,q)$ up to conjugation (cf.~\cite{Su}). It follows that all
the subfield subgroups of the same order in $\PSL(2,q)$ lie in at
most two conjugacy classes, whereas all the subfield subgroups of
the same order in $\PGL(2,q)$ are conjugate. Similarly, all the
subgroups isomorphic to either $A_4$, $S_4$ or $A_5$ lie in at most
two conjugacy classes in $\PSL(2,q)$, and they belong to the same
conjugacy class in $\PGL(2,q)$.

Let $G=\PSL(2,q)$ (or $\PGL(2,q)$), and let $\tilde G = \PGL(2,q)$,
so that $G \leq \tilde G$. Let $H$ be a subfield subgroup or a small
subgroup of $G$, and let $\tilde H$ denote the normalizer $N_{\tilde
G}(H)$.

The following proposition summarizes some useful properties
regarding conjugation of elements and subgroups in $G$. These
properties quickly follow from the aforementioned classification of
elements and subgroups up to conjugation.

\begin{prop} \label{conj}
Let $G$,$\tilde G$,$H$ and $\tilde H$ be as above. Then the
following hold.
\begin{enumerate}
\item
If $K \leq G$ and $K \cong H$ then there exists $\tilde g \in \tilde
G$ such that $K = H^{\tilde g}$. Moreover, in this case also $\tilde
K = \tilde H^{\tilde g}$.
\item
If $w,w' \in H$ are conjugate in $\tilde G$, then $w$ and $w'$ are
already conjugate in $\tilde H$.
\item
If $K \geq H$ then $K$ is either a subfield subgroup or a small
subgroup. Moreover, $\tilde K \geq \tilde H$.
\end{enumerate}
\end{prop}

%%%%%%%%%%%%%%%%%%%%%%%%%%%%%%%%%%%%%%%%%%%%%%%%%%%%%%%

\subsection{Spread} \label{spread}

Recall that a $2$-generated group $G$ has \emph{spread $m$} if for
every subset of $m$ nontrivial group elements $g_1,\ldots,g_m$,
there exists an element $h \in G$ such that $\langle g_1,h \rangle =
\ldots = \langle g_m,h \rangle = G$. We say that $G$ has \emph{exact
spread $m$} if it has spread $m$ but not $m+1$. Also recall that a
generating $k$-tuple $(g_1,\dots,g_k)$ for $G$ is called
\emph{redundant} if there is some $i$ such that the group generated
by $\{g_j\}_{j \neq i}$ is equal to $G$.

\begin{prop} \label{spread}
Let $G$ be a group of spread $2$ and let $k\geq 3$ be an integer.
Then all redundant generating $k$-tuples of $G$ belong to the same
connected component in $\tilde \Gamma_k(G)$.
\end{prop}
\begin{proof} See~\cite[Lemma 2.8]{E}.
\end{proof}

Recently it has been proved that any finite simple group $G$ has
spread $2$~\cite{BGK,GK}. In particular, the results in~\cite{BGK}
imply that the groups $\PSL(2,q)$ and $\PGL(2,q)$ (for $q>3$) have
spread $2$. Moreover, the exact spread of $\PSL(2,q)$ is computed
in~\cite{BW} for $q \neq 7$.

\begin{prop} \label{spread_PSL}
Suppose $G=\PSL(2,q)$, where $q \neq 7$. Then the exact spread of
$G$, which we denote by $\alpha$, is given as follows:
\[
\begin{array} {|c|c c c c c c|}
\hline q & 2 & 3 & 5 & 9 & \geq 4 \text{ even} & \geq 11 \text{ odd} \\
\hline \alpha & 3 & 4 & 2 & 2 & q-2 & \left\{
\begin{matrix} q-1 & q \equiv 1 \pmod 4 \\ q-4 & q \equiv 3 \pmod 4
\end{matrix}\right.\\
\hline
\end{array}
\]
\end{prop}

The exact spread of $\PSL(2,7)$ is known to be at least $3$
(see~\cite{BW}), although determining the precise value remains an
open problem. Similar methods can be used to compute the exact
spread of $\PGL(2,q)$ (see~\cite{thesis}).

\begin{prop} \label{spread_PGL}
Suppose $G=\PGL(2,q)$, where $q$ is an odd prime power. Then the
exact spread of $G$, which we denote by $\alpha$, is given as
follows:
\[
\begin{array}{|c|c c c c|}
\hline q & 3 & 5 & 7 & \geq 9\\
\hline \alpha & 1 & 2 & 4 & q-4 \\
\hline
\end{array}
\]
\end{prop}

%%%%%%%%%%%%%%%%%%%%%%%%%%%%%%%%%%%%%%%%%%%%%%%%%%%%%%%

\section{Proof of theorem \ref{main}}

Let $G=\PSL(2,q)$ or $\PGL(2,q)$, where $q=p^e$ and $p$ is an odd
prime. Since $\PSL(2,q)$ and $\PGL(2,q)$ both have spread $2$ it is
enough to connect any generating $k$-tuple to a redundant one (by
Proposition~\ref{spread}).

For clarity of presentation, Lemmas~\ref{normalizer}--~\ref{eq_ord}
below are stated and proved for the minimal size of a generating
$k$-tuple (i.e. $k=3$ or $4$ as appropriate). However, they are
valid for any larger value of $k$.

The proof of Theorem~\ref{main} involves the adaptation of some of
the techniques of Evans~\cite{E}. For the convenience of the reader,
we review these techniques in detail.

We note that unfortunately we could not achieve the desired result
with $k=3$. The main obstacle lies in Step 3 below, where techniques
other than those in \cite{E} had to be used.

%%%%%%%%%%%%%%%%%%%%%%%%%%%%%%%%%%%%%%%%%%%%%%%%%%%%%%%

\subsection{Step 1: Finding an element of order different than $2$ or $p$}
In the first step we connect a given generating tuple $(w,y,z)$ (or
$(w,x,y,z)$) to a tuple $(w',y',z')$ (or $(w',x',y',z')$) where $w'$
is a split or a non-split element, of order different than $2$.

\begin{lemma} \label{normalizer}
\emph{(Adaptation of \cite[Lemma 4.2]{E})}. Let $(w,y,z) \in V_3(G)$
where $w \neq 1$. Then $(w,y,z)$ is connected to $(w,y',z')$, where
$y',z' \notin N_G (\langle w \rangle)$. Moreover, we may assume that
$y'$ and $z'$ are not of order $2$.
\end{lemma}
\begin{proof}
Note that $y$ and $z$ cannot both be in the normalizer $N_G (\langle
w \rangle)$, since $G$ does not normalize $\langle w \rangle$. If $y
\in N_G (\langle w \rangle)$ and $z \notin N_G (\langle w \rangle)$,
then $yz \notin N_G (\langle w \rangle)$. Thus, we can connect
$(w,y,z) \rightarrow (w,yz,z)$ and $yz,z \notin N_G (\langle w
\rangle)$. Therefore, we may now assume that in the generating
$3$-tuple $(w,y,z)$ we have $y,z \notin N_G (\langle w \rangle)$.

If $y$ is of order $2$, then $wy$ is of order different than $2$.
Indeed, if $wywy = 1$, then $w^{-1} = ywy = w^{y}$, so $w^y \in
\langle w \rangle$, implying that $y$ normalizes $\langle w
\rangle$, a contradiction. Therefore, if $y$ is of order $2$, we can
connect $(w,y,z) \rightarrow (w,wy,z)$, where $|wy| \neq 2$ and $wy
\notin N_G (\langle w \rangle)$. We can apply the same argument for
$z$ if necessary.
\end{proof}

\begin{lemma} \label{not_2_p}
If $x,y \in G$ are two non-commuting elements of order $p$, then
either there exists some $i$ with $|xy^i| \neq 2,p$, or $p=3$ and
$\langle x,y \rangle \cong A_4$.
\end{lemma}
\begin{proof}
Since $x$ and $y$ are unipotent, they are both stabilizers of points
in $\mathbb{P}^1(q)$. Assume that $x$ fixes $a$ and $y$ fixes $b$,
then $a$ and $b$ are distinct points because $x$ and $y$ do not
commute. Therefore, we may change coordinates, and assume that $a,b$
are the images of the vectors $[1,0], [0,1] \in \mathbb{F}_q^2$. In
these coordinates, $x$ and $y$ are the images of the matrices
\[
    X =
    \left( \begin{matrix}
    1 & \lambda \\ 0 & 1
    \end{matrix} \right)
    \quad , \quad
    Y =
    \left( \begin{matrix}
    1 & 0 \\ \mu & 1
    \end{matrix} \right)
    \qquad
    (\lambda,\mu \in \mathbb{F}_q),
\]
under the natural projection map.

Therefore,
\[
    XY^j =
    \left( \begin{matrix}
    1+j \lambda\mu & \lambda \\ j\mu & 1
    \end{matrix} \right).
\]
This induces an element $xy^j \in G$ of order $p$ if and only if it
has exactly one eigenvalue. However, the characteristic polynomial
of $XY^j$ is
\[
    p(t) = (t-1-j\lambda\mu)(t-1)-j\lambda\mu = t^2 -
    (2+j\lambda\mu)t + 1
\]
and its discriminant is $j\lambda\mu(4+j\lambda\mu)$. Therefore,
$XY^j$ is unipotent if and only if $j = 0$ or $4+j\lambda\mu \equiv
0 \pmod p$.

If $p>3$ then we can always find an integer $0<j<p-1$ for which
$|xy^j|, |xy^{j+1}| \neq p$. If either $|xy^j|$ or $|xy^{j+1}|$ does
not equal $2$, then we can take $i=j$ or $i=j+1$ accordingly, and we
are done. Otherwise, $x \in N_G{(\langle y \rangle)}$ and this
implies that $x$ and $y$ commute, a contradiction.

If $p=3$, the argument above shows that $xy$ and $xy^2$ cannot both
be of order $2$. We conclude that either there exists some $i$ with
$|xy^i| \notin \{2,3\}$, or we are in one of the following two
cases:
\begin{enumerate}
\item $1= x^3 = y^3 = (xy)^2 = (xy^2)^3$, and then $\langle x,y \rangle
\cong A_4$; or
\item $1= x^3 = (y^2)^3 = (xy^2)^2 = (xy)^3$, and then
$\langle x,y^2 \rangle = \langle x,y \rangle \cong A_4$.
\end{enumerate}

This completes the proof of the lemma.
\end{proof}

\begin{lemma} \label{2_p}
Assume that $p>3$ and let $(w,y,z) \in V_3(G)$. Then $(w,y,z)$ is
connected to $(w',y',z')$, where $|w'| \neq 2,p$.
\end{lemma}
\begin{proof}
If no such $(w',y',z')$ exists then the orders $|w|, |y|$ and $|z|$
all equal $2$ or $p$.

\emph{Case (i):} $|w| = |y| = |z| =2$.

If $|wy| = |yz| = |wz| = 2$, then $w,y$ and $z$ all commute with
each other, a contradiction since $G$ is non-abelian. Therefore, we
may assume that $|wy| \neq 2$ and connect $(w,y,z) \rightarrow
(wy,y,z)$. If $|wy| \neq p$ we are done, otherwise $|wy| = p$ and we
are in the situation of Case (ii).

\emph{Case (ii):} At least one of $|w|, |y|, |z|$ equals $p$.

Without loss of generality we may assume that $|w| = p$. By Lemma
\ref{normalizer}, we may also assume that $|y|,|z| \neq 2$. If $|y|
\neq p$ or $|z| \neq p$ we are done. Otherwise, $|w| = |y| = |z| =
p$, and we are in the situation of Case (iii).

\emph{Case (iii):} $|w| = |y| = |z| =p$.

It is not possible that $w,y$ and $z$ all commute with each other.
Thus we may assume that $w$ and $y$ do not commute. By Lemma
\ref{not_2_p}, since $p>3$, there exists some $i$ such that $|wy^i|
\neq 2,p$. Thus we may connect $(w,y,z) \rightarrow (wy^i, y, z)$
and we are done.
\end{proof}

\begin{lemma} \label{2_3}
Assume that $p=3$, and let $(w,x,y,z) \in V_4(G)$. Then $(w,x,y,z)$
is connected to $(w',x',y',z')$, where $|w'| \neq 2,3$.
\end{lemma}
\begin{proof}
We proceed as in the proof of Lemma~\ref{2_p}. For a contradiction,
suppose that $|w|, |x|, |y|$ and $|z|$ all equal $2$ or $3$.

\emph{Case (i):} At least one of $|w|,|x|,|y|,|z|$ equals $3$.

Without loss of generality we may assume that $|w|=3$ and that
$x,y,z \notin N_G(\langle w \rangle)$, by Lemma \ref{normalizer}. In
particular, none of the elements $x$,$y$ or $z$ commute with $w$.
The proof of Lemma \ref{not_2_p} reveals that either $|\alpha w|
\neq 3$ or $|\alpha w^2| \neq 3$ for each $\alpha \in \{x,y,z\}$. If
there exists some $i$ with $|\alpha w^i| \neq 2,3$ then we are done.
Otherwise, we may connect $(w,x,y,z) \rightarrow (w,xw^i,yw^j,zw^k)
= (w,x',y',z')$, where $|x'|=|y'|=|z'|=2$, and so we reduce to the
situation of Case (ii).

\emph{Case (ii):} $|x|=|y|=|z|=2$.

The subgroup $\langle x,y \rangle$ is dihedral. More precisely, if
$t=|xy|$ then $\langle x,y \rangle$ is a dihedral subgroup of order
$2t$. By inspecting the list of subgroups presented in
\S\ref{subgroups}, we see that $t$ divides $3^m \pm 1$ for some $m$,
therefore $t \neq 3$. If $t \neq 2$, then we are done, otherwise, if
$t=2$, then $x$ and $y$ commute. Similarly, the groups $\langle
x,z\rangle$ and $\langle y,z \rangle$ are dihedral. Therefore, we
reduce to the case where $\langle x,y,z \rangle$ is elementary
abelian of order $8$, and this is not one of the optional subgroups
in \S\ref{subgroups}.
\end{proof}

%%%%%%%%%%%%%%%%%%%%%%%%%%%%%%%%%%%%%%%%%%%%%%%%%%%%%%%

\subsection{Step 2: Eliminating the structural subgroups}
In the second step we connect a generating $3$-tuple $(w,y,z)$ to
$(w,y'z')$, where $w$ is a split or a non-split element of order
different than $2$, and $\langle w,w^{y'} \rangle$ and $\langle
w,w^{z'} \rangle$ are not structural subgroups.

\begin{lemma} \label{struct_split_1} \emph{(Adaptation of \cite[Lemma 4.6]{E})}.
Let $(w,y,z) \in V_3(G)$ where $w$ is a split element of order
different than $2$. Suppose that $(w,y,z)$ is not connected to a
redundant $3$-tuple. Then $(w,y,z)$ is connected to $(w,y',z')$,
where $\langle w,y' \rangle$ and $\langle w,z' \rangle$ are not
structural subgroups, and neither $y'$ nor $z'$ are of order $2$.
\end{lemma}

\begin{proof}
By Lemma \ref{normalizer} we may assume that $y,z \notin N_G
(\langle w \rangle)$ and $|y|,|z| \neq 2$. Since $w$ is a split
element it fixes two points on $\mathbb{P}_1(q)$, say $a$ and $b$.
Suppose that $\langle w,y \rangle$ is a structural subgroup. The
proper structural subgroups of $G$ that contain $w$ are $N_G
(\langle w \rangle)$, $G_a$ and $G_b$, together with certain
subgroups of these. Since $y$ does not normalize $w$, we see that
$\langle w,y \rangle$ is a subgroup of $G_a$ or $G_b$. Without loss
of generality we may assume that $\langle w,y \rangle \leq G_a$. We
distinguish two cases.

\emph{Case (i): $y$ fixes $a$ only.} Recall that $G_a$ is a Borel
subgroup, and let $P_a$ denote the $p$-Sylow part of $G_a$, i.e.
$G_a = N_G(P_a)$. Then all the elements in $P_a$ are of order $p$
and so they have only one fixed point, namely $a$.

Now, $y \in P_a$ and $w \notin P_a$, thus $yw \notin P_a$, so that
$yw$ fixes an additional point, $d$ say, and clearly $yw \notin N_G
(\langle w \rangle)$. If $yw$ is not of order $2$, the
transformation $(w,y,z) \rightarrow (w,yw,z)$ gives a $3$-tuple of
the sort to be considered in Case (ii) below. Otherwise if $yw$ is
of order $2$, then $yw^2$ is not of order $2$, and again, $y \in
P_a$ and $w^2 \notin P_a$. Therefore $yw^2 \notin P_a$, so $yw^2$
fixes an additional point, $d$ say, and clearly $yw^2 \notin N_G
(\langle w \rangle)$. Thus the transformation $(w,y,z) \rightarrow
(w,yw^2,z)$ gives a $3$-tuple of the type we now consider in Case
(ii).

\emph{Case (ii): $y$ fixes an additional point $d \neq a$, i.e. $y
\in G_{a,d}$.} If $d=b$ then $\langle w,y \rangle \in G_{a,b}$ and
we note that $G_{a,b}$ is cyclic. Thus $\langle w,y \rangle$ is
cyclic and $(w,y,z)$ is connected to a redundant tuple by Dunwoody's
result~\cite{Du2} on solvable groups. Therefore we may assume that
$b \neq d$.

Now, $za \neq a$, otherwise $\langle w,y,z \rangle \leq G_a$, a
contradiction. Also, if $za=b$ and $z^{-1}a=b$, then $w^za=a$ and
$w^zb=b$, thus $w^z \in G_{a,b}$ which is cyclic. Since $w$ and
$w^z$ are of the same order, we have $w^z \in \langle w \rangle$,
hence $z \in N_G (\langle w \rangle)$, a contradiction.

Therefore we can define $z'=z^{\pm 1}$ such that $z'a \neq a$ and
$z'a \neq b$. Let $g_i = z'w^i$ for $i=0,1,2$, and observe that
$y^{g_i^{-1}} \in G_{g_i a, g_i d} = G_{z'a, g_i d}$. Suppose that
$g_i d = g_j d$ where $0 \leq i < j \leq 2$. Then $w^{j-i}d=d$ and
thus $j=i$ since $d \neq a,b$. It follows that $g_k d \neq a,b$ for
some $k=0,1,2$. Now, $y^{g_k^{-1}} \in G_{z'a, g_k d}$, and so
$y^{g_k^{-1}}$ fixes $z'a$ and $g_kd$ but no other point. Hence,
$y^{g_k^{-1}}$ has two fixed points, neither of which is $a$ or $b$.
Therefore, by the list of subgroups in \S\ref{subgroups}, $\langle
w,y^{g_k^{-1}} \rangle$ is not a structural subgroup. Clearly we can
transform $(w,y,z) \rightarrow (w,y,z') \rightarrow (w,y,g_k)
\rightarrow (w,y^{g_k^{-1}},g_k) \rightarrow (w,y^{g_k^{-1}},z')
\rightarrow (w,y^{g_k^{-1}},z)$. Let $y'=y^{g_k^{-1}}$ and observe
that $y'$ is again of order different than $2$.

The above argument shows that if $\langle w,y \rangle$ is a
structural subgroup then we can transform $(w,y,z) \rightarrow
(w,y',z)$, where $\langle w,y' \rangle$ is not a structural subgroup
and $|y'| \neq 2$. Similarly, if $\langle w,z \rangle$ is a
structural subgroup then we can repeat the same argument for
$(w,y',z)$ and obtain the desired result, using the fact that now
$y'a \neq a$, $y'b \neq b$ and $y' \notin N_G (\langle w \rangle)$.
\end{proof}

\begin{lemma} \label{struct_split_2} \emph{(Adaptation of \cite[Lemma 4.7]{E})}.
Let $(w,y,z) \in V_3(G)$ where $w$ is a split element of order
different than $2$. Suppose that $(w,y,z)$ is not connected to a
redundant $3$-tuple. Then $(w,y,z)$ is connected to $(w,y',z')$,
where $\langle w,w^{y'} \rangle$ and $\langle w,w^{z'} \rangle$ are
not structural subgroups.
\end{lemma}

\begin{proof}
By Lemma \ref{struct_split_1} we may assume that $\langle w,y
\rangle$ and $\langle w,z \rangle$ are not structural subgroups and
that $|y|,|z| \neq 2$. In particular, we may assume that $y,z \notin
N_G (\langle w \rangle)$. Since $w$ is split, there exists $a,b \in
\mathbb{P}_1(q)$ such that $w \in G_{a,b}$. Then $w^y \in
G_{y^{-1}a,y^{-1}b}$, so if $\langle w,w^y \rangle$ is a structural
subgroup, then it is a subgroup of $G_a$ or $G_b$. Without loss of
generality, we may assume that $\langle w,w^y \rangle \leq G_a$.

However, $y^{-1} \notin G_a$ as $\langle w,y \rangle$ is not a
structural subgroup. Since $w^y \in G_a$ it follows that
$y^{-1}b=a$. Now $z \notin N_G (\langle w \rangle)$, so either $za
\neq b$ or $z^{-1}a \neq b$. Set $z'=z^{\pm 1}$ such that $z'a \neq
b$. Note that since $\langle w,z' \rangle$ is not a structural
subgroup, it has no fixed points, thus $z'a \neq a$ and $z'b \neq
b$.

Define $g_i = z'w^iy^{-1}$ for $i=0,1,2$, and suppose that $g_i a=
g_j a$ for some $0 \leq i < j \leq 2$. Then $y^{-1}a =
w^{j-i}y^{-1}a$, so $w^{j-i}$ fixes $y^{-1}a$. Hence either
$y^{-1}a=a$ or $y^{-1}a=b$, since these are the only points fixed by
$w$. However $y^{-1}a \neq a$ and $y^{-1}b=a$, thus $y^{-1}a=b$
implies that $w^ya=a$ and $w^yb=b$, thus $w^y \in G_{a,b}$ which is
cyclic. Since $w$ and $w^y$ are of the same order, $w^y \in \langle
w \rangle$, implying that $y \in N_G (\langle w \rangle)$, a
contradiction.

Therefore $g_0 a, g_1 a, g_2 a$ are distinct, and so $g_k a \neq
a,b$ for some $k=0,1,2$. Now $w^{g_k^{-1}} \in G_{g_k a, g_k b} =
G_{g_k a, z'a}$, and so $w^{g_k^{-1}}$ fixes $g_k a$ and $z'a$ but
no other point. Hence, $w^{g_k^{-1}}$ has two fixed points, neither
of which is $a$ or $b$. We deduce that $\langle w,w^{g_k^{-1}}
\rangle$ is not a structural subgroup. Clearly we can transform
$(w,y,z) \rightarrow (w,y,z') \rightarrow (w,g_k^{-1},z')
\rightarrow (w,g_k^{-1},z) = (w,y',z)$, and $\langle w,w^{y'}
\rangle$ is not a structural subgroup. If $\langle w, w^z \rangle$
is a structural subgroup, we can repeat the same argument for
$(w,y',z)$ and obtain the desired result.
\end{proof}

\begin{lemma} \label{struct_non-split} \emph{(Adaptation of \cite[Lemma 4.8]{E})}.
Let $(w,y,z) \in V_3(G)$ where $w$ is a non-split element of order
different than $2$. Suppose that $(w,y,z)$ is not connected to a
redundant $3$-tuple. Then $(w,y,z)$ is connected to $(w,y',z')$,
where $\langle w,w^{y'} \rangle$ and $\langle w,w^{z'} \rangle$ are
not structural subgroups.
\end{lemma}

\begin{proof}
Since $w$ is a non-split element of order different than $2$, the
only structural subgroups of $G$ that contain $w$ are the subgroups
of $N_G (\langle w \rangle)$. By Lemma \ref{normalizer} we may
assume that $y,z \notin N_G (\langle w \rangle)$. Since $w$ is of
order different than $2$, if $w^y \in N_G (\langle w \rangle)$,
which is a dihedral group, then $w^y \in \langle w \rangle$.
Therefore $y \in N_G (\langle w \rangle)$, a contradiction.
Consequently, $\langle w,w^y \rangle$ is not a structural subgroup.
Similarly, we deduce that $\langle w,w^z \rangle$ is also
non-structural.
\end{proof}

%%%%%%%%%%%%%%%%%%%%%%%%%%%%%%%%%%%%%%%%%%%%%%%%%%%%%%%

\subsection{Step 3: Connecting to a redundant tuple}
In this final step we connect a generating $4$-tuple $(w,x,y,z)$,
for which $\langle w,w^x \rangle$, $\langle w,w^y \rangle$ and
$\langle w,w^z \rangle$ are not structural subgroups, to a redundant
$4$-tuple.

By Lemmas \ref{struct_split_2} and \ref{struct_non-split}, we may
assume that $L_1 = \langle w,w^{x} \rangle$, $L_2 = \langle w,w^{y}
\rangle$, and $L_3 = \langle w,w^{z} \rangle$ are not structural
subgroups. Let $K_1 = \langle w,x \rangle$, $K_2 = \langle w,y
\rangle$, and $K_3 = \langle w,z \rangle$; and let $H_1 = \langle
w,x,y \rangle$, $H_2 = \langle w,x,z \rangle$, and $H_3 = \langle
w,y,z \rangle$. Note that $K_1$,$K_2$ and $K_3$ are not structural
subgroups, since they contain a non-structural subgroup. Similarly,
$H_1$, $H_2$ and $H_3$ are also non-structural.

If one of the $H_i$ is isomorphic to a small subgroup (i.e. $A_4$,
$S_4$ or $A_5$) then we are done since any generating $3$-tuple of
such a group is connected to a redundant one. For $A_4$ and $S_4$,
this follows from~\cite{Du2}, while the result for $A_5 \cong
\PSL(2,5)$ was obtained by Gilman~\cite{Gi}.

Therefore, we may assume that $H_1$, $H_2$ and $H_3$ are subfield
subgroups, and denote $\tilde H_1 \cong \PGL(2,q_1)$, $\tilde H_2
\cong \PGL(2,q_2)$ and $\tilde H_3 \cong \PGL(2,q_3)$ (where
$q_1,q_2$ and $q_3$ are odd prime powers that divide $q$). Without
loss of generality, we may also assume that $q_1 \leq q_2 \leq q_3$.

\begin{lemma} \label{diff_ord} \emph{(Adaptation of Case (ii) of \cite[Lemma 4.9]{E})}.
With respect to the above notation, if $q_1 < q_3$ then one of the
following holds:
\begin{enumerate}
\renewcommand{\theenumi}{\roman{enumi}}
\item $(w,x,y,z)$ is connected to a redundant tuple;
\item $(w,x,y,z)$ is connected to $(w,x',y,z)$, where $M_1 =
\langle w,x',y \rangle$, $M_2 = \langle w,x',z \rangle$ and $M_3 =
\langle w,y,z \rangle = H_3$ satisfy the bound
\[
    |\tilde M_1| + |\tilde M_2| +|\tilde M_3| >
    |\tilde H_1| + |\tilde H_2| + |\tilde H_3|
\]
with $\tilde M_i = N_{\PGL(2,q)}(M_i)$.
\end{enumerate}
\end{lemma}

\begin{proof}
Since $w,w^x \in L_1$, there exist some $\tilde d \in \tilde L_1$
such that $w^x = w^{\tilde d}$. Set $\tilde u = x \tilde d^{-1}$, so
$\tilde u \in C_{\tilde K_1} (w) $. Let $C_{H_3}(w) = \langle c
\rangle$ and observe that $\langle c,\tilde u \rangle$ is cyclic
since it is a subgroup of $C_{\tilde G} (w)$. Thus, by \cite[Lemma
4.3]{E} there exists an integer $n$ such that $\langle c,\tilde u
\rangle = \langle c^n \tilde u \rangle$. Now, $c \in H_3 = \langle
w, y, z \rangle$, so we may connect $(w,x,y,z) \rightarrow (w,c^n x,
y,z)$.

Define $M_1 = \langle w, c^n x , y \rangle$ and $M_2 = \langle w,
c^n x , z \rangle$, and note that $M_1$ and $M_2$ are non-structural
since they both contain $\langle w, w^{c^n x} \rangle = \langle w,
w^x \rangle = L_1$. Since $\tilde d \in L_1$ we also have $\tilde d
\in \tilde M_1 \cap \tilde M_2$, therefore $\tilde M_1$ and $\tilde
M_2$ contain $(c^n x) \tilde d^{-1} = (c^n x) (x^{-1} \tilde u) =
c^n \tilde u$. Further, since $c \in \langle c^n \tilde u \rangle$,
we deduce that $c \in \tilde M_1 \cap \tilde M_2$, hence $C_{H_3}(w)
= \langle c \rangle \leq C_{\tilde M_i}(w)$ for $i=1,2$.

If $M_1$ or $M_2$ is isomorphic either to $A_4 \cong \PSL(2,3)$,
$S_4 \cong \PGL(2,3)$ or $A_5 \cong \PSL(2,5)$, then we can get a
redundant tuple and (i) holds. Otherwise, $\tilde M_i =
\PGL(2,q_i')$ where $q'_1,q'_2>3$ are odd prime powers that divide
$q$, and $|C_{\tilde M_i}(w)| = q_i' \pm 1$ whereas
$|C_{H_3}(w)|=q_3 \pm 1$ or $\frac{1}{2}(q_3 \pm 1)$. An analysis of
the inequality
\[
\frac{1}{2}(q_3 - 1) \leq |C_{H_3}(w)| \leq |C_{\tilde M_1}(w)| \leq
q'_1+1
\]
shows that if $q_3 > q'_1$ then $q_3=9$ and $q_1'=3$, a
contradiction. Therefore, $q_3 \leq q'_1$ and similarly $q_3 \leq
q'_2$. The hypotheses $q_1 \leq q_2 \leq q_3$ and $q_1 < q_3$ now
imply that $|\tilde M_1| + |\tilde M_2| + |\tilde H_3|
> |\tilde H_1| + |\tilde H_2| + |\tilde H_3|$, so that (ii) holds.
\end{proof}

\begin{lemma} \label{eq_ord}
With respect to the above notation, if $q_1=q_3$ then $(w,x,y,z)$ is
connected to a redundant tuple.
\end{lemma}
\begin{proof}
In this case, $|\tilde H_1| = |\tilde H_2| = |\tilde H_3| =
|\PGL(2,q_1)|$. At least two of $|H_1|$, $|H_2|$ and $|H_3|$ are
equal (either to $|\PSL(2,q_1)|$ or to $|\PGL(2,q_1)|$), so without
loss of generality, we may assume that $|H_1|$ = $|H_2|$.

Next, we note that $H_1 \cap H_2 \geq \langle w, x \rangle = K_1$.
Since $H_1$ and $H_2$ are subfield subgroups of the same order,
Proposition \ref{conj} implies that there exists $\tilde g \in
\tilde G$ such that $H_1 ^ {\tilde g} = H_2$. Since $K_1$ is a
subgroup of $H_1$, it follows that ${K_1}^{\tilde g}$ is a subgroup
of $H_2$. However, $K_1$ is also a subgroup of $H_2$. Therefore,
$K_1$ and ${K_1}^{\tilde g}$ are conjugate in $\tilde H_2$. Thus
there exists $\tilde h \in \tilde H_2$ such that ${K_1}^{\tilde g
\tilde h} = K_1$. Therefore $\tilde g \tilde h \in N_{\tilde G}
(K_1) = \tilde K_1 \leq \tilde H_2$, and thus $\tilde g \in \tilde
H_2 = N_{\tilde G} (H_2)$. Now, $H_1 = H_2 ^{\tilde g ^{-1}} = H_2$
and $\langle H_1,H_2 \rangle = G$, so one of the generators $y$ or
$z$ is redundant.
\end{proof}

\begin{corr}
Any generating $4$-tuple $(w,x,y,z) \in V_4(G)$ is connected to a
redundant tuple.
\end{corr}
\begin{proof}
By Lemmas \ref{2_p} and \ref{2_3} we may assume that $w$ is either a
split or non-split element of order different than $2$. Therefore,
by Lemmas \ref{struct_split_2} and \ref{struct_non-split}, we can
connect any generating $4$-tuple $(w,x,y,z)$ to a tuple
$(w,x',y',z')$, such that $\langle w,w^{x'} \rangle$, $\langle
w,w^{y'} \rangle$ and $\langle w,w^{z'} \rangle$ are not structural
subgroups. Let $H_1 = \langle w,x',y' \rangle$, $H_2 = \langle
w,x',z' \rangle$, and $H_3 = \langle w,y',z' \rangle$, and note that
$H_1$, $H_2$ and $H_3$ are not structural subgroups. Moreover, we
may assume that the $H_i$ are subfield subgroups.

Among all tuples that are connected to $(w,x',y',z')$ we can take
the tuple for which $|\tilde H_1| + |\tilde H_2| + |\tilde H_3|$ is
maximal. However, if $|\tilde H_1|,|\tilde H_2|$ and $|\tilde H_3|$
are not all equal to each other, then Lemma \ref{diff_ord} yields a
contradiction, and if $|\tilde H_1| = |\tilde H_2| = |\tilde H_3|$,
then Lemma \ref{eq_ord} yields the desired connectivity to a
redundant tuple.
\end{proof}

This completes the proof of Theorem~\ref{main}.

%%%%%%%%%%%%%%%%%%%%%%%%%%%%%%%%%%%%%%%%%%%%%%%%%%%%%%%

\vspace{\bigskipamount}

\end{document}